\documentclass[12pt, reqno]{amsart}

\usepackage{epsfig,psfrag,subfigure}
\usepackage{epstopdf}
\usepackage{amsmath,amssymb,amsfonts,latexsym, mathrsfs,euscript}
\usepackage{graphicx,graphics,cite}
\usepackage{wrapfig}
\usepackage[usenames]{color}
\usepackage{geometry}                % See geometry.pdf to learn the layout options. There are lots.
\usepackage{colortbl,booktabs,ctable,multirow}
\usepackage{epstopdf}
\usepackage{subfig}

\DeclareGraphicsExtensions{.pdf,.png,.jpg}

\makeatletter \@addtoreset{figure}{section}
\def\thefigure{\thesection.\@arabic\c@figure} \def\fps@figure{h, t}
\@addtoreset{table}{section}
\def\thetable{\thesection.\@arabic\c@table}\def\fps@table{h, t}
\@addtoreset{equation}{section}

\makeatother

\newcommand\bs{\boldsymbol}
\newcommand\beq{\begin{equation}}
\newcommand\eeq{\end{equation}}

\title[A numerical scheme for double-gyre model]{A conservation formulation and a numerical algorithm for the double-gyre nonlinear shallow-water model}

\author{Dongyang Kuang}
\address{Department of Mathematics, University of Wyoming, Laramie, WY 82071, USA}
\email{dkuang@uwyo.edu}

\author{ Long Lee}
\address{Department of Mathematics, University of Wyoming, Laramie, WY 82071, USA}
\email{llee@uwyo.edu}

\date{\today}

\begin{document}

\maketitle
%\markboth{Dongyang Kuang and Long Lee}{A conservation formulation and a numerical algorithm for the double-gyre nonlinear shallow-water model}
%\title{A conservation formulation and a numerical algorithm for the double-gyre nonlinear shallow-water model}
%
%\author[Dongyang Kuang and Long Lee]{Dongyang Kuang and Long Lee\corrauth}
%\address{Department of Mathematics, University of Wyoming, Laramie, WY 82071, USA.}
%\email{{\tt llee@uwyo.edu} (Long Lee)}

%%%%% Begin Abstract %%%%%%%%%%%
\begin{abstract}
We present a conservation formulation and a numerical algorithm for the reduced-gravity shallow-water equations on a beta plane, subjected to a constant wind forcing that leads to the formation of double-gyre circulation in a closed ocean basin. The novelty of the paper is that we reformulate the governing equations into a nonlinear hyperbolic conservation law plus source terms. A second-order fractional-step algorithm is used to solve the reformulated equations. In the first step of the fractional-step algorithm, we solve the homogeneous hyperbolic shallow-water equations by the wave-propagation finite volume method. The resulting intermediate solution is then used as the initial condition for the initial-boundary value problem in the second step. As a result, the proposed method is not sensitive to the choice of viscosity and gives high-resolution results for coarse grids, as long as the Rossby deformation radius is resolved. We discuss the  boundary conditions in each step, when no-slip boundary conditions are imposed to the problem. We validate the algorithm by a periodic flow on an f-plane with exact solutions. The order-of-accuracy for the proposed algorithm is tested numerically. We illustrate a quasi-steady-state solution of the double-gyre model via the height anomaly and the contour of stream function for the formation of double-gyre circulation in a closed basin. Our calculations are highly consistent with the results reported in the literature. Finally, we present an application, in which the double-gyre model is coupled with the advection equation for modeling transport of a pollutant in a closed ocean basin. 
\end{abstract}
%%%%% end %%%%%%%%%%%

\begin{description}
\item[{\footnotesize\bf keywords:} ]
{\footnotesize double-gyre, reduced-gravity shallow-water equations, wave-propagation algorithm, fractional-step algorithm.}
\end{description}

%%%% Start %%%%%%
\section{Introduction}

The two-dimensional shallow-water equations govern the fluid motion in a thin layer. They can be used as a rational approximation to the three-dimensional Euler equations, with the assumption of hydrostaticity and shallow water depth (compared with the horizontal length scale).  When wind forcing and latitude-dependent Coriolis forces are included, these equations represent a simple model for describing the depth-average dynamics of the oceans. Furthermore, if we include a Laplacian diffusion in the equations and impose Dirichlet boundary conditions on the velocity field, in particular the no-slip conditions, the equations are often used to simulate a mid-latitude closed ocean basin. In this paper we focus on a reduced-gravity shallow-water model formulated for studying the behavior of western boundary currents (WBCs) in mid latitudes \cite{JJG95}. In this ocean model water is assumed to consist of two layers of fluid, a single active layer of fluid of constant density $\rho$ and variable thickness $h(x, y, t)$, overlying a deep and motionless layer of density $\rho+\Delta\rho$. Consequently, the motion of the upper layer represents the gravest baroclinic mode \cite{JJG95}.  The model equations in non-conservation form are
\begin{equation}\label{eq:double_gyre}
\begin{split}
&\frac{\partial h}{\partial t}+\frac{\partial(uh)}{\partial x} +\frac{\partial(vh)}{\partial y}=0,\\
&\frac{\partial u}{\partial t} +u\frac{\partial u}{\partial x} + v\frac{\partial u}{\partial y}=-g_r\frac{\partial h}{\partial x} + (f_0+\beta y)v + \nu\nabla ^2 u +F^{u},\\ 
&\frac{\partial v}{\partial t} +u\frac{\partial v}{\partial x} + v\frac{\partial v}{\partial y}=-g_r\frac{\partial h}{\partial y} - (f_0+\beta y)u + \nu\nabla ^2 v +F^{v},\\ 
\end{split}
\end{equation}
where $(u, v)$ is the velocity filed, $h$ is the height field, $g_r=(\Delta\rho / \rho)g$ is the reduced gravity, and $g$ is the acceleration of gravity. $(F^{u}, F^{v})$ is the external forcing term, such as the wind forcing  \cite{JJG95, JPM97, Pedlosky, PJM96}. With the imposition of no-slip boundary conditions on the velocity field (the height field is allowed to assume any value on the boundaries), equations (\ref{eq:double_gyre}) describe a wind-driven, closed basin on a $\beta$ plane. The equations are normally referred to as the double-gyre, wind-driven shallow-water model. This model is a convenient test bed for studying mid-latitude ocean dynamics \cite{JJG95, JPM97, PJM96 }. 

The numerical algorithm MPDATA (Multidimensional Positive Definite Advection Transport Algorithm) has long been used to solve geophysical flows, such as flow governed by Eq. (\ref{eq:double_gyre}). MPDATA is a two-pass scheme that preserves positive definite scalar transport functions with small oscillations \cite{mpdata84, mpdata93, mpdata98}.  Technically, the method belongs to the same class of non-oscillatory Lax-Wendroff algorithms such as FCT \cite{FCT}, TVD\cite{TVD}, and ENO \cite{ENO}. Nevertheless, MPDATA was primarily developed for meteorological applications. The method focuses on reducing the implicit viscosity of the donor cell scheme, while retaining the virtues of positivity, low phase error, and simplicity of upstream differencing. However, the disadvantage of MPDATA is that the basic MPDATA is too diffusive, and enhanced MPDATA is too expensive \cite{mpdata98}. We compare a basic MPDATA implementation described in \cite{hayder06} with the proposed algorithm for the double-gyre model in Section \ref{sec3}. For a thorough review of MPDATA, we refer the readers to \cite{mpdata98}.

Aimed at improving the resolution and accuracy, a type of multi-scale finite difference method was developed in \cite{JPM97, PJM96 } for solving equations (\ref{eq:double_gyre}). The multi-scale method, or enslaved finite-difference method makes use of properties of the governing equations in the absence of time derivatives to reduce the overall truncation errors without changing the order of spatial discretization, nor the time step restriction of the time integrator. This means that the enslaved scheme effectively increases the spatial resolution of the given algorithm without changing its temporal stability or memory requirements. However, the enslaved scheme could be sensitive to the viscosity values used in the calculation of solving the shallow-water double-gyre model for some time integrators. Especially for numerical approximations with resolution near the Rossby deformation radius. For example, it is reported in \cite{JPM97} that for Rossby deformation radius $\approx 52-75$ km, the implementation of an enslaved scheme using the leapfrog time integrator could be numerically unstable for explicit viscosity values less than $\nu =1000$ $\text{m}^2\,\text{s}^{-1}$ for the resolution $\Delta x =40$ km, and  $\nu =750$ $\text{m}^2\,\text{s}^{-1}$ for the resolution $\Delta x =20$ km.  For solving the double-gyre model it is common for this class of schemes that to maintain numerical stability, the viscosity needs to be increased as the grid resolution is decreased. \cite{JPM97}.

%Other parameters used in their calculations are the same as those used in our numerical experiments, but the viscosity used in our calculation can be fixed as small as $\nu =200$ $\text{m}^2\,\text{s}^{-1}$ for these grid sizes.  

In this paper, we propose a stable method for solving the double-gyre model. We rewrite the governing equations into a conservation form with source terms. A fractional-step algorithm is used to solves the new formulation. In the first step, the hyperbolic equations are solved by the high-resolution wave-propagation method developed by LeVeque \cite{bib:clawpack}. Then the resulting intermediate values are used as the initial conditions for the initial-boundary problem. 
The fractional-step strategy has proven to be efficient and stable for solving the Navier-Stokes equations and other fluid models \cite{bib:ll03, bib:ll10}. 

We organize the rest of the paper as follows. In Section \ref{sec2}, we present the conservation form of the double-gyre shallow-water model. Then we introduce a fractional-step method to solve the equations and discuss the boundary conditions in each step. In Section \ref{sec3}, we verify the algorithm by an exact solution of a period flow on an f-plane. We  show that numerically the method is second-order accurate. Then we use the algorithm to study an upper-ocean double-gyre circulation in a closed ocean basin. We compute the height anomaly for the formation of double-gyre circulation, and compare the results with the literature values computed by the enslaved finite-difference schemes \cite{PJM96} and the traditional methods of backward Euler and  centered finite-difference \cite{JJG95}. The results are highly consistent. Finally, we present an example, in which the double-gyre model is coupled with the advection equation for modeling transport of a pollutant in a closed ocean basin. This example demonstrates the flexibility of the proposed method to couple with other equations that require high-resolution results for the monitored quantity, such as a passive tracer in fluid. 

\section{The fractional-step algorithm}\label{sec2}

The model equations (\ref{eq:double_gyre}) can be written in their conservation form
\begin{equation}\label{eq:double_gyre_conv}
\begin{split}
&\frac{\partial h}{\partial t}+\frac{\partial(uh)}{\partial x} +\frac{\partial(vh)}{\partial y}=0,\\
&\frac{\partial (hu)}{\partial t} +\frac{\partial}{\partial x}\left(hu^2+\frac{1}{2}g_rh^2\right)+\frac{\partial}{\partial y}\left(huv\right) = (f_0+\beta y) hv + h (\nu\nabla^2 u) +h F^{u},\\
&\frac{\partial (hv)}{\partial t} +\frac{\partial}{\partial x}\left(huv\right)+\frac{\partial}{\partial y}\left(hv^2+\frac{1}{2}g_rh^2\right) =- (f_0+\beta y) hu + h (\nu\nabla^2 v) +h F^{v},
\end{split}
\end{equation}
where $hu$ and $hv$ are the momenta in $x$ and $y$ directions, and 
\begin{equation}\label{eq:p}
P(h)= \frac{1}{2}g_rh^2
\end{equation}
is the hydrostatic equation of state with a reduced gravity. Equations (\ref{eq:double_gyre_conv}) represent a system of two-dimensional hyperbolic conservation law with a source term,
\begin{equation}\label{eq:hcl_2d}
q_t+f(q)_x+g(q)_y=\psi(q, \tilde{q}),
\end{equation}
where
\begin{equation}
\begin{split}
&q=\begin{bmatrix}
h \\ hu \\ hv
\end{bmatrix}, \qquad
f(q)=\begin{bmatrix}
hu \\ hu^2+\frac{1}{2}g_rh^2 \\ huv
\end{bmatrix}, \qquad
g(q)=\begin{bmatrix}
hv \\ huv \\ hv^2+\frac{1}{2}g_rh^2
\end{bmatrix},\\
&\tilde{q}=
\begin{bmatrix}
h\\u\\v
\end{bmatrix}, \qquad
\psi(q, \tilde{q})=\begin{bmatrix}
0\\ (f_0+\beta y)hv + h(\nu\nabla^2 u) +h F^{u}\\- (f_0+\beta y)hu + h(\nu\nabla^2 v) +h F^{v}
\end{bmatrix}.
\end{split}
\end{equation}

At first glance, equations (\ref{eq:double_gyre_conv}) seem to be inconsistent in the treatment of the stress tensor parametrization. The natural variable for the momentum equations is $q$, so in principle the assumed eddy viscosity parametrization should also be expressed in terms of $q$, instead of $u$ and $v$. However, scaling the non-conservation ``advective'' form of the shallow-water equations (\ref{eq:double_gyre}) leads to the geostrophic balance between the horizontal velocity and the horizontal pressure gradient, i.e. the gradient of height field. In other words, the principal geostrophic balance is between the Coriolis force and the height (pressure) gradient, not the dissipative term\cite{Pedlosky}. The conservation formulation (\ref{eq:double_gyre_conv}) preserves the principal geostrophic balance, and the balance is now in the form of the momentum variable $q$.

We propose a fractional-step method, also known as operator splitting, for Eq. (\ref{eq:hcl_2d}) that simply alternates solving the following two problems:
\begin{equation}\label{eq:two_problems}
\begin{split}
&\text{Problem A:}\quad q_t + f(q)_x + g(q)_y = 0;\\
&\text{Problem B:}\quad q_t  =  \psi(q, \tilde{q}).\\
\end{split}
\end{equation}
Problem A is a homogeneous  conservation law that can be solved by the high-resolution finite volume method developed in \cite{bib:clawpack}.  After spatial discretization, Problem B is treated as a simple system of ordinary differential equations (ODEs) that can be solved by a standard time integrator. Since $h_t = 0$ in Problem B, we can further simplify Problem B by letting 
\begin{equation}
q_1=\begin{bmatrix}
hu \\ hv
\end{bmatrix},\quad
\tilde{q}_1=\begin{bmatrix}
u \\ v 
\end{bmatrix},
\end{equation}
and Problem B becomes
\begin{equation}\label{eq:Prob_B_matrix}
\frac{\partial q_1}{\partial t}= R q_1 + S (\tilde{q}_1, h),
\end{equation}
where $R$ is a $2\times 2$ constant matrix and $S$ is a vector function of $\tilde{q}_1$ and $h$. The forms of $R$ and $S$ are explicitly written in Eq. (\ref{eq:step2}).

If both Problem A and B are solved over one time step $\Delta t$, this is the so-called Godunov splitting for a fractional-step method. The splitting error of the Godunov splitting is $O(\Delta t)$ in theory. In practice, however, the error is smaller than $O(\Delta t)$ \cite{CL02}. The Strang splitting is a slight modification of the Godunov splitting and yields second-order accuracy generally \cite{bib:clawpack}. The difference between the Godunov splitting and the Strang splitting is that the Strang splitting starts and ends with a half time step $\Delta t /2$ on Problem A. In between the first and the last time steps, the Strang splitting is the same as the Godunov splitting. That is, Problem B and A are solved alternately over one time step $\Delta t$.  The splitting error of Strang splitting is $O(\Delta t^2)$. To be more specific, basically for the Godunov splitting we solve the two sub-problems sequentially, like (A) $\longrightarrow$ (B),  by using the time increments $\{\Delta t ,\,\Delta t\}$ in each time step, respectively, and for the Strang splitting in each time step we solve the two sub-problems in a sequence of (A) $\longrightarrow$ (B) $\longrightarrow$  (A), by using the time increments $\{\frac{\Delta t}{2},\,\Delta t,\,\frac{\Delta t}{2}\}$. After combining the the cycles, $\{\frac{\Delta t}{2},\,\Delta t,\,\frac{\Delta t}{2}\}$, $\{\frac{\Delta t}{2},\,\Delta t,\,\frac{\Delta t}{2}\}$,\dots, $\{\frac{\Delta t}{2},\,\Delta t,\,\frac{\Delta t}{2}\}$, the Strang splitting is the same as the Godunov splitting, except the Strang splitting uses $\frac{\Delta t}{2}$ for solving Problem A in the very beginning, as well as the very end. Moreover, Yoshida \cite{YH}  introduced a systematic method to construct arbitrary even-order time accurate splitting schemes. The Strang splitting is a modification of the fist member of the Yoshida's method.

%Note that if $f_0$ and $\beta$ are both small in ${\psi}(q, \tilde{q})$, Problem B is a stiff system of ODEs,  a method that is not $L$-stable, such as the trapezoidal method, should not be used to solve the problem \cite{bib:LeVeque_diff}. A stiff ODE solver will be required for the problem.  
 
%\subsection{The algorithm}
 
Let the computational domain be $\Omega$ and the boundary of the domain be $\partial\Omega$. Let $C$ be a two-dimensional grid cell $\Delta x \times \Delta y$ and  $q$ be the solution the partial differential equations. Let $Q_{i,j}^{n}$ be an approximation to the cell average of $q$ over the cell $C_{i,j}$ at time $t=t^{n}$, i.e. 
\begin{equation}
Q_{i,j}^{n}=\frac{1}{\Delta x\Delta y}\int_{C_{i,j}} q(x, y, t^{n})dx dy.
\end{equation} 
The cell averaged value is placed at the cell center.
Suppose that the boundary conditions for the velocity field $u$ and $v$ are prescribed, the fractional-step method is described as follows:
\begin{itemize}
\item {\bf Step 1}:  Given $Q^{n}$, the semi-discrete system of equations arising from Problem A has the form
\begin{equation}\label{eq:step1}
\begin{split}
&\frac{Q^{m}-Q^{n}}{\Delta t}  + F(Q^{m}, Q^{n}) +  G(Q^{m}, Q^{n}) = 0,\\
&Q^{m}\,\,\text{on}\,\,{\partial\Omega},\,\,\text{the boundary conditions are given}.
\end{split}
\end{equation}
Solve the above system by the wave-prorogation finite volume method to obtain $Q^{m}$. We briefly describe the multidimensional wave-prorogation finite volume method as follows. Problem A, the hyperbolic shallow-water equations, can be written as a quasi-linear equations
\beq
q_t +f'(q)q_x +g'(q)q_y=0,
\eeq
where 
\beq
f'(q) = A(h,u,v) =\left( \begin{array}{ccc}
0 & 1 & 0 \\
-u^2+g_rh& 2u & 0 \\
-uv& v & u 
\end{array} \right),
\eeq
and
\beq
g'(q) = B(h,u,v) =\left( \begin{array}{ccc}
0 & 0 & 1 \\
-uv & v & u \\
-v^2 + g_rh& 0 & 2v
\end{array} \right).
\eeq
Let $c=\sqrt{g_rh}$ be the speed of gravity wave. The matrix $A$ has eigenvalues and eigenvectors
\beq\label{eq:e-value-A}
\begin{split}
& \lambda^{x_1} =u-c,\quad \lambda^{x_2} =u,\quad \lambda^{x_3} =u+c \\
& r^{x_1}=\left[ \begin{array}{c}
1\\
u-c\\
v
\end{array} \right],\quad
r^{x_2}=\left[ \begin{array}{c}
0\\
0\\
1
\end{array} \right],\quad
r^{x_3}=\left[ \begin{array}{c}
1\\
u+c\\
v
\end{array} \right],
\end{split}
\eeq
while the matrix $B$ 
has eigenvalues and eigenvectors
\beq\label{eq:e-value-B}
\begin{split}
& \lambda^{x_1} =v-c,\quad \lambda^{x_2} =v,\quad \lambda^{x_3} =v+c \\
& r^{x_1}=\left[ \begin{array}{c}
1\\
u\\
v-c
\end{array} \right],\quad
r^{x_2}=\left[ \begin{array}{c}
0\\
-1\\
0
\end{array} \right],\quad
r^{x_3}=\left[ \begin{array}{c}
1\\
u\\
v+c
\end{array} \right].
\end{split}
\eeq
For the wave-propagation algorithm, the updating formula over a time step $\Delta t$ is
\beq\label{eq:update_Q}
\begin{split}
Q^{m}_{i, j} = Q^{n}_{i, j} &- \frac{\Delta t}{\Delta x}\left(\mathcal{A}^{+}\Delta Q^{n}_{i-1/2,j}+ \mathcal{A}^{-}\Delta Q^{n}_{i+1/2,j} \right)\\
& - \frac{\Delta t}{\Delta y}\left(\mathcal{B}^{+}\Delta Q^{n}_{i,j-1/2}+ \mathcal{B}^{-}\Delta Q^{n}_{i,j+1/2} \right)\\
& - \frac{\Delta t}{\Delta x}\left(\tilde{F}_{i+1/2, j} - \tilde{F}_{i-1/2, j} \right) - \frac{\Delta t}{\Delta y}\left(\tilde{G}_{i, j+1/2} - \tilde{G}_{i, j-1/2} \right). 
\end{split}
\eeq
The second and the third terms on the right-hand-side of Eq. (\ref{eq:update_Q}) are the fluctuations, while the fourth and the fifth terms are the correction fluxes. Both fluctuations and correction fluxes are computed by using the approximate Riemann solver (or the Roe solver) that averags the waves and speeds (corresponding to the eigenvectors and eigenvalues in Eqs. (\ref{eq:e-value-A}) \& (\ref{eq:e-value-B})) by the Roe average. Detailed information about the actual representations of the fluctuations and correction fluxes can be found in \cite{bib:clawpack}, pp 471--474.

\vskip 1cm 

\item {\bf Step 2}: From {\bf Step 1}, we obtain
\begin{equation}
Q^{m}=\begin{bmatrix}
H^{m} \\ (HU)^{m} \\ (HV)^{m},
\end{bmatrix}, \qquad
H^{n+1}=H^{m}.
\end{equation}

The semi-discretized equations for (\ref{eq:Prob_B_matrix}), arising by using the centered-difference scheme for the spatial derivatives, has the form
\begin{equation}\label{eq:step2}
\begin{split}
\frac{\partial (HU)_{i, j}}{\partial t}&=(f_0+\beta y)(HV)_{i,j} + \\ &\nu H_{i,j}\left(\frac{U_{i-1, j} - 2 U_{i, j} + U_{i+1, j} }{(\Delta x)^2} + \frac{U_{i, j-1} - 2 U_{i, j} + U_{i, j+1} }{(\Delta y)^2} \right)+ H_{i,j} F^{u}_{i.j},\\
\frac{\partial (HV)_{i, j}}{\partial t} &= (f_0+\beta y)(HU)_{i,j} +\\ & \nu H_{i,j}\left(\frac{V_{i-1, j} - 2 V_{i, j} + V_{i+1, j} }{(\Delta x)^2} + \frac{V_{i, j-1} - 2 V_{i, j} + V_{i, j+1} }{(\Delta y)^2} \right)+ H_{i,j} F^{v}_{i.j},
\end{split}
\end{equation}
where 
\beq\label{eq:U_and_V}
U_{i,j}=\frac{(HU)_{i,j}}{H_{i,j}},\quad V_{i,j}=\frac{(HV)_{i,j}}{H_{i,j}},\quad\text{for}\,\, H_{i,j}\ne 0,\,\,i,\,j=1\cdots N.
\eeq
If $H_{i,j}=0$, it means that the water depth is zero, which is not physically meaningful. Equation (\ref{eq:step2}) is a system of ODEs with $2N$ dimensions, where $N$ is the number of cells used in (\ref{eq:step1}). The initial conditions are $(HU)_{i.j}=(HU)^m_{i,j}$ and $(HV)_{i.j}=(HV)^m_{i,j}$ and the final time is $t^{n+1}=t^{n}+\Delta t$. The prescribed boundary conditions for the velocity field are employed in this step. 
Note that with a sufficiently small time step, an explicit  $p$-stage, $p^{th}$-order Runge-Kutta method, $p > 1$, is a $A$-stable method for  solving equation (\ref{eq:step2}), for which $U_{i,j}=U^{m}_{i,j}$, $V_{i, j}=V^{m}_{i,j}$, and $H_{i, j}=H^{m}_{i,j}$ \cite{bib:LeVeque_diff}.  

\end{itemize}
It is worth noting that in order to simulate a closed ocean basin, no-slip boundary conditions ($u=0$ and $v=0$) are usually prescribed for the non-conservation model equation (\ref{eq:double_gyre}). For the fractional-step algorithm, however,  we require two sets of boundary conditions: one for $Q^{m}$ in the first step and one for $U$, $V$ in the second step. Naturally, the no-slip boundary conditions are employed in the second step (\ref{eq:step2}), while the choice of the boundary conditions for $Q^{m}$ in the first step must reflect the physical interpretation of no-slip boundary conditions. We choose solid-wall boundary conditions for $Q^{m}$. The key observation of a solid wall is that at the boundary $x=a$,
\begin{equation}\label{eq:solid_wall_u}
u(a, y, t)= 0,\quad hu(a,y,t)=0.
\end{equation}
Similarly, a solid wall at the boundary $y=b$ is
\begin{equation}\label{eq:solid_wall_v}
v(x, b, t)= 0,\quad hv(x,b,t)=0.
\end{equation}
To achieve the solid-wall conditions (\ref{eq:solid_wall_u}) and  (\ref{eq:solid_wall_v}), in each time step the ghost-cell values in the second-order finite volume wave-propagation algorithm are set to be
\begin{equation}\label{eq:ghost_cells}
\begin{split}
\text{For}\,\, Q^{m}_0&:\quad H^{m}_0 =H^{m}_{1},\quad (HU)^{m}_0=-(HU)^{m}_1, \quad (HV)^{m}_0=-(HV)^{m}_1\\
\text{For}\,\, Q^{m}_{-1}&:\quad H^{m}_{-1} =H^{m}_{2},\quad (HU)^{m}_{-1}=-(HU)^{m}_{2}, \quad (HV)^{m}_{-1}=-(HV)^{m}_{2}.
\end{split}
\end{equation}
Formula (\ref{eq:ghost_cells}) imposes a necessary symmetry for achieving the solid-wall conditions (\ref{eq:solid_wall_u}) and (\ref{eq:solid_wall_v}) \cite{bib:clawpack}. 
We remark that because we discretize the Laplacian by the five-point centered-difference scheme, and both $U$ and $V$ are computed at cell centers, we are not using the no-slip conditions $V=0$ at the vertical walls and $U=0$ at the horizontal walls. Instead, we use the ghost cell values of $U$ and $V$, which are computed based on Eq. (\ref{eq:U_and_V}) and Eq (\ref{eq:ghost_cells}). Our choice of the ghost cell values enforces the boundary condition prescribed for the non-conservation equations. We also remark that the homogeneous hyperbolic shallow-water equations (\ref{eq:hcl_2d}) is solved by the high-resolution wave propagation algorithms developed by LeVeque \cite{bib:clawpack} in this study. The algorithms can easily be replaced by other efficient anti-diffusion shock-capturing schemes, such as the algorithms developed in \cite{shu1, shu2} and many others, for which we do not attempt to provide a detailed list.  We hope to emphasize that this study focuses on introducing a new formulation for the double-gyre shallow-water model and a numerical implementation for solving the formulation. It is not our intention to develop a new efficient anti-diffusion shock-capturing scheme, neither to develop a new algorithm for solving the hyperbolic conservation laws with source terms. We demonstrate that as a result of combining the new formulation and the fractional-step algorithm, we obtain a stable method that is not sensitive to the kinematic viscosity and the grid refinement for the double-gyre shallow-water model.

\section{Numerical investigation}\label{sec3}

\subsection{Periodic flows on an f-plane}

We validate the proposed algorithm by examining a periodic flow on a constant f-plane (i.e., the Coriolis force does not depend on latitude and thus $\beta=0$). We introduce the following dimensionless variables:
\begin{equation}\label{eq:variables}
u^{*}=\frac{u}{U},\,\,v^{*}=\frac{v}{U},\,\,h^{*}=\frac{h}{H_0},\,\,x^{*}=\frac{x}{L},\,\,y^{*}=\frac{y}{L},\,\,t^{*}=\frac{t}{L/U},
\end{equation}
where $U$ is the scale of velocity, $L$ is the typical length scale, and $H_0$ is the scale of water height. Substituting the above dimensionless variables into equations (\ref{eq:double_gyre_conv}) results in the following scaled system of equations for the double-gyre shallow-water model (we drop `*' herein and after) :
\begin{equation}\label{eq:non_dimen}
\begin{split}
&\frac{\partial h}{\partial t}+\frac{\partial(uh)}{\partial x} +\frac{\partial(vh)}{\partial y}=0,\\
&\frac{\partial (hu)}{\partial t} +\frac{\partial}{\partial x}\left(hu^2+\frac{1}{2}F_{r}^{-2} h^2\right)+\frac{\partial}{\partial y}\left(huv\right) = \frac{1}{R_{0}} hv + \frac{1}{Re}h \nabla^2 u +h F^{u},\\
&\frac{\partial (hv)}{\partial t} +\frac{\partial}{\partial x}\left(huv\right)+\frac{\partial}{\partial y}\left(hv^2+\frac{1}{2}F_{r}^{-2} h^2\right) =- \frac{1}{R_{0}} hu + \frac{1}{Re}h \nabla^2 v +h F^{v},
\end{split}
\end{equation}
%$$\text{where}\,\,F_0=\frac{U}{\sqrt{gH_0}}\,\,\text{is the Froude number},\,\, Re=\frac{LU}{\nu}\,\, \text{is the Reynolds number},\,\, \text{and}\,\,R_0=\frac{U}{Lf}\,\,\text{is the Rossby number}.$$
%where $F_{r}=\displaystyle\frac{U}{\sqrt{g_rH_0}}$ is the Froude number, $Re=\displaystyle\frac{LU}{\nu}$ is the Reynolds number, and $R_0=\displaystyle\frac{U}{Lf}$ is the the Rossby number.
where $F_{r}=U / \sqrt{g_rH_0}$ is the Froude number, $Re=LU / \nu$ is the Reynolds number, and $R_0=U / Lf$ is the Rossby number.
Consider the solution ansatz 
\begin{equation}\label{eq:ansatz}
\begin{split}
u(x,y,t)&=(\eta +\epsilon\sin(\omega t))\cos(2\pi x)\sin(2\pi y),\\
v(x,y,t)&=-(\eta +\epsilon\sin(\omega t))\sin(2 \pi x)\cos(2\pi y),\\
h(x,y)&=\exp(\cos(2 \pi x)\cos(2\pi y)),
\end{split}
\end{equation}
where the parameters $\eta$, $\epsilon$, and $\omega$ control the contribution of spatial and temporal derivatives in the solution.  Substituting the ansatz into equation (\ref{eq:non_dimen}), we obtain the forcing terms $F^u$,
\begin{equation}\label{eq:fu}
\begin{split}
F^{u}= & \epsilon\omega\cos(\omega t)\cos(2\pi x)\sin(2\pi y)-2\pi(\eta+\epsilon\sin(\omega t))^2\sin(2\pi x)\cos(2\pi x)\\ & + \frac{8\pi^2}{Re}(\eta+\epsilon\sin(\omega t))\cos(2\pi x)\sin(2\pi y) + \frac{1}{R_0}(\eta+\epsilon\sin(\omega t))\sin(2\pi x)\cos(2\pi y)\\ & -\frac{2\pi}{Fr^2}\sin(2\pi x)\cos(2\pi y)\exp(\cos(2\pi x)\cos(2\pi y)),
\end{split}
\end{equation}
and $F^v$,
\begin{equation}\label{eq:fv}
\begin{split}
F^{v}= & - \epsilon\omega\cos(\omega t)\sin(2\pi x)\cos(2\pi y)-2\pi(\eta+\epsilon\sin(\omega t))^2\sin(2\pi y)\cos(2\pi y)\\ & - \frac{8\pi^2}{Re}(\eta+\epsilon\sin(\omega t))\sin(2\pi x)\cos(2\pi y) + \frac{1}{R_0}(\eta+\epsilon\sin(\omega t))\cos(2\pi x)\sin(2\pi y)\\ & -\frac{2\pi}{Fr^2}\cos(2\pi x)\sin(2\pi y)\exp(\cos(2\pi x)\cos(2\pi y)).
\end{split}
\end{equation}
Note that the solution ansatz is independent of the dimensionless parameters $F_r$, $R_0$, and $Re$. In principle, for this test problem, the solution behavior of the proposed fractional-step method should be insensitive to the choice of these parameters, if the following conditions are satisfied: 1)  the CFL condition in the first step of solving the hyperbolic equation and 2) the stability restriction of the 2-stage, second-order Runge-Kutta method used to solve equation (\ref{eq:step2}).
%if the CLF condition in the first step of solving the hyperbolic equation and the stability restriction of the 2-stage, second-order Runge-Kutta method used to solve the second step, equation (\ref{eq:step2}), are both satisfied. 
We choose $Re=100$, $R_0=0.1$, and $Fr=2$ for our simulations.  For the required boundary conditions in (\ref{eq:step1}) and (\ref{eq:step2}), we impose periodic boundary conditions in both steps for periodic flow.

\begin{figure}[h]
\begin{center}
(a)\includegraphics[width=2.55in]{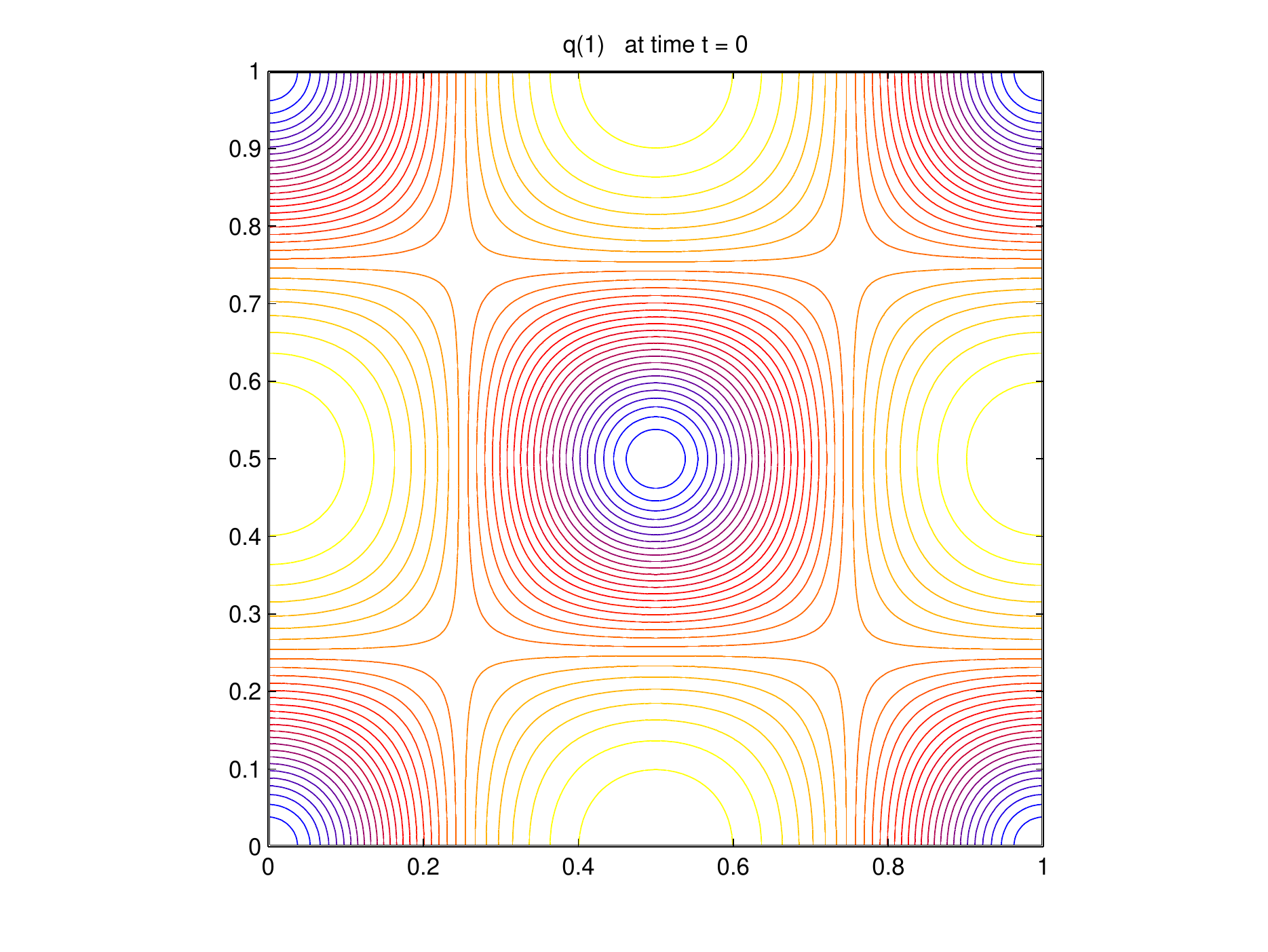}
(b)\includegraphics[width=2.55in]{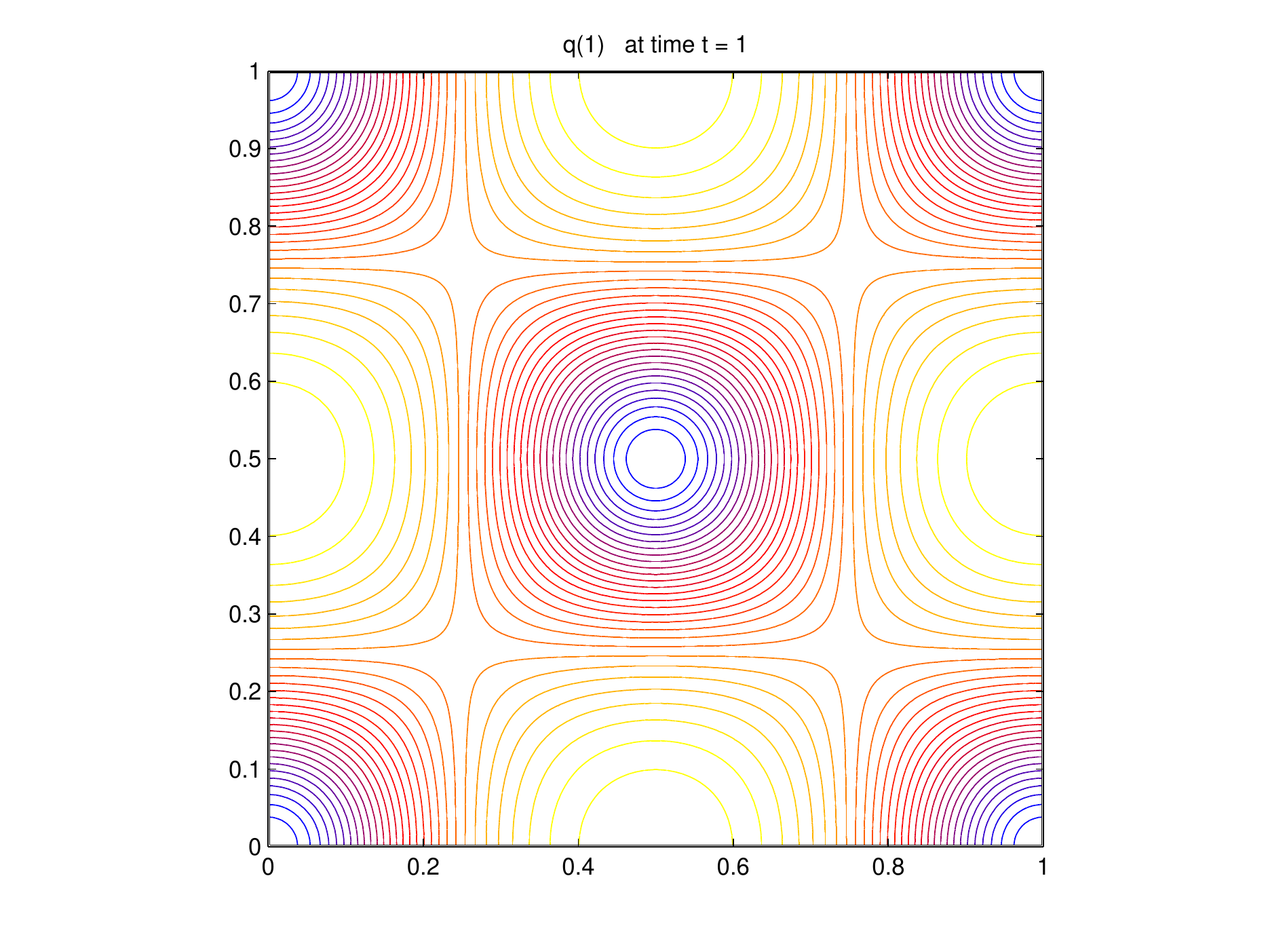}
\end{center}
\caption{The height field of a periodic flow on an f plane. (a) Exact solution of the height field. 30 contour lines are used for values between 0.36794 and 2.7179. (b) Numerical solution at $t=1$ with $250\times 250$ cells. 30 contour lines are used for values between 0.36796 and 2.7183.}
\label{fig:f1}
\end{figure}

For computational domain $[ 0,1]\times[0,1]$,  Table \ref{tab1} shows the grid refinement study for the fractional-step method. We compute the error of the height field between the exact solution and the numerical solution at the final time $t=1$, with the (finite) $l_2$-norm 

\begin{equation}
\label{eq:error}
||\epsilon||= \sqrt{\frac{1}{N^2}\sum_{i=1}^{N }\sum_{j=1}^{N}\epsilon_{i,j}^2},
\end{equation}
where $N$ is the number of grid cells in one direction. 
Since the error goes down by four (on average) when we refine the grid, it provides evidence that the proposed method is second-order accurate. We use the parameters $\eta=0.1$, $\epsilon=0.9$, and $\omega=\pi / 20$. Note that the time step $\Delta t$ in this calculation is chosen so that when we refine the grid in both $x$ and $y$ directions, the time step used for the fine gird is $1/4$ of that used for the coarse grid. We start with $\Delta t =0.025$ for $N=10$. The Strang-splitting method is used for the calculation. No limiters are used in the first step. Figure \ref{fig:f1}(a) is the exact solution of the height field. 30 contour lines are used for values between 0.36794 and 2.7179. Figure \ref{fig:f1}(b) is the numerical solution at $t=1$ with $250\times 250$ cells. 30 contour lines are used for values between 0.36796 and 2.7183.

\begin{table}[htbp]
\begin{centering}%\footnotesize
\caption{Convergence rate for the fractional-step algorithm.}
\label{tab1}
\begin{tabular}{ccccccc} 
\hline
$N$ & 10& 20 & 40 & 80 & 160 & 250\\
\hline
$||H_N-h_{exact}||$& 5.133e-2 &1.203e-2& 3.304e-3 & 8.718e-4 & 2.300e-4  & 9.467e-5\\
\hline
order & &2.09 &1.87  & 1.92 & 1.92 & 1.99\\
\hline
\end{tabular}
\end{centering}
\end{table}

Unlike the enslaved schemes developed in  \cite{PJM96, JPM97}, the proposed numerical method is independent of the magnitude of the time dependent contribution to the solution. That is, the accuracy of numerical solution and the efficiency of the algorithm are independent of the choice of $\eta$ and $\epsilon$. Table \ref{tab:t2} shows errors of the computed solutions in the $l_2$-norm for the horizontal velocity $u$, and the elapsed CPU times for various choices of $\eta$ and $\epsilon$. As expected, the numerical experiments show that the solution behavior of the proposed algorithm is insensitive to the choice of $\eta$ and $\epsilon$.  Note that the absolute error increases as $\eta$ increases, due to the fact that  the magnitude of $u$ increases as $\eta$ increases. The numerical experiments use a $50\times 50$ grid, while the parameter $\omega = \pi/10$ and the final run time is $t=5$. %The computation was carried out using MacBook Pro with 2.16 GHz Intel Core Duo Processor. 

%In \cite{PJM96}, the authors indicate that if $\eta + \epsilon =1$ and $\epsilon$ is large ( e.g. $\epsilon =0.9$ in this simulation), then the enslaved scheme is computationally less efficient than the original scheme for the leapfrog time integrator. 

\begin{table}[htbp]
\begin{centering}%\footnotesize
\caption{Errors of computed solutions for $u$ and the elapsed CPU times for various choices of $\eta$ and $\epsilon$ .}
\label{tab:t2}
\begin{tabular}{cccccc} 
\hline
$\eta+\epsilon = 1$ & $\eta = 0.1$ & $\eta=0.3$ & $\eta=0.5$ &$\eta=0.7$ & $\eta=0.9$\\
\hline
$||U_N-u_{exact}||$& 4.07e-3 &4.19e-3& 4.36e-3 & 4.54e-3 & 4.70e-3   \\
\hline
CPU time (sec) &35.39 &35.30 & 35.33  & 35.30 & 35.30 \\
\hline
\end{tabular}
\end{centering}
\end{table}

\subsection{Upper-ocean double-gyre model}\label{sec:DG}

To demonstrate the strength of the method that combines the new formulation and the fractional-step algorithm, we examine the geophysical flow that describes a closed basin flow on a $\beta$-plane, subjected to zonal winds. With reduced gravity, the model resembles a two-layer ocean basin whose upper layer is driven by a zonal wind stress \cite{JJG95}, e.g. the external forcing term in equation (\ref{eq:double_gyre}) is the imposed wind forcing given by the curl of the wind stress,
\begin{equation}\label{eq:wind_forcing}
\begin{split}
F^{u}& = -\frac{\tau_0}{\rho H_0}\cos\left(\frac{2\pi y}{L}\right),\\
F^{v}& = 0,
\end{split}
\end{equation} 
where $\tau_0$ is the wind stress, $\rho$ is the water density, $L$ is the domain length in the North-South direction, and $H_0$ is the initial upper-layer depth. The parameter values used in the simulations are listed in Table \ref{tab:t3}. These values are chosen to closely match of those in \cite{JJG95, PJM96} for comparison. 
\begin{table}[htbp]
\caption{Model parameters.}
\label{tab:t3}
\begin{centering}\footnotesize
\begin{tabular}{lcl} 
\hline
Coriolis parameter & & $f_0=5.0\times 10^{-5}$s$^{-1}$\\
$f=f_0+\beta y$ & & $\beta=1.875\times 10^{-11}$\\
Wind stress & & $\tau_0 = 0.11$ N m$^{-2}$\\
Kinematic viscosity & & $\nu=300$ m$^{2}$s$^{-1}$\\
Upper-layer density & & $\rho=1000$ kg m$^{-3}$\\
Reduced gravity & & $g_r=0.03$ ms$^{-2}$\\
Initial upper layer depth & & $H_0 = 500$ m\\
Domain & & $D=1000$ km (East-West)\\
& & $L=2000$ km (North-South)\\
\hline
\end{tabular}
\end{centering}
\end{table}

\begin{figure}[htpb]
\begin{center}
(a)\includegraphics[width=1.38in]{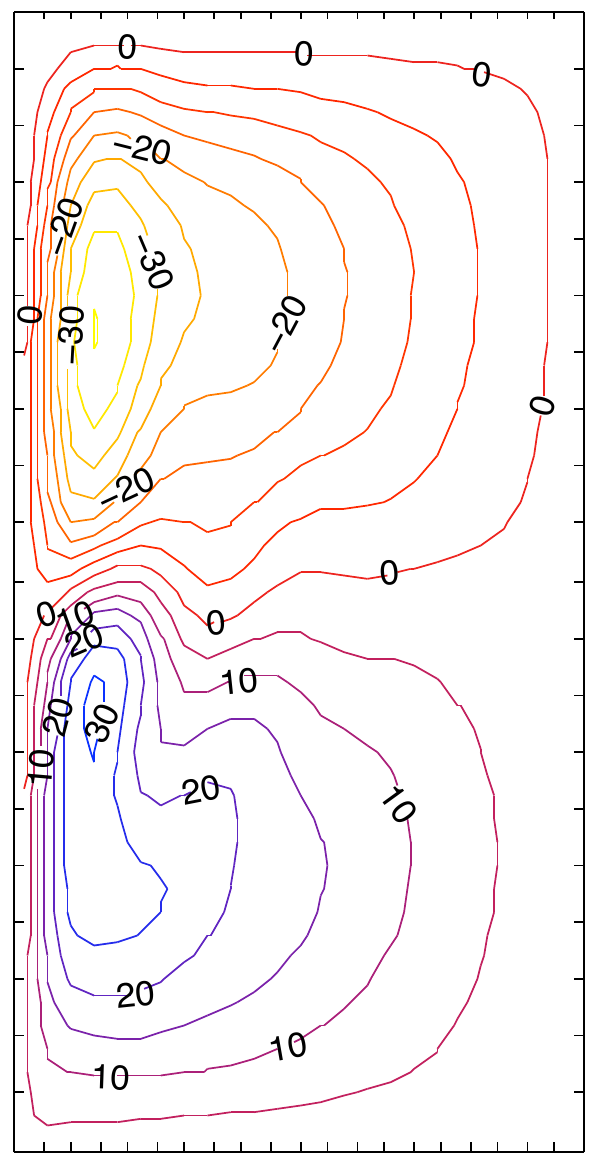}\hspace{7mm}
(b)\includegraphics[width=1.4in]{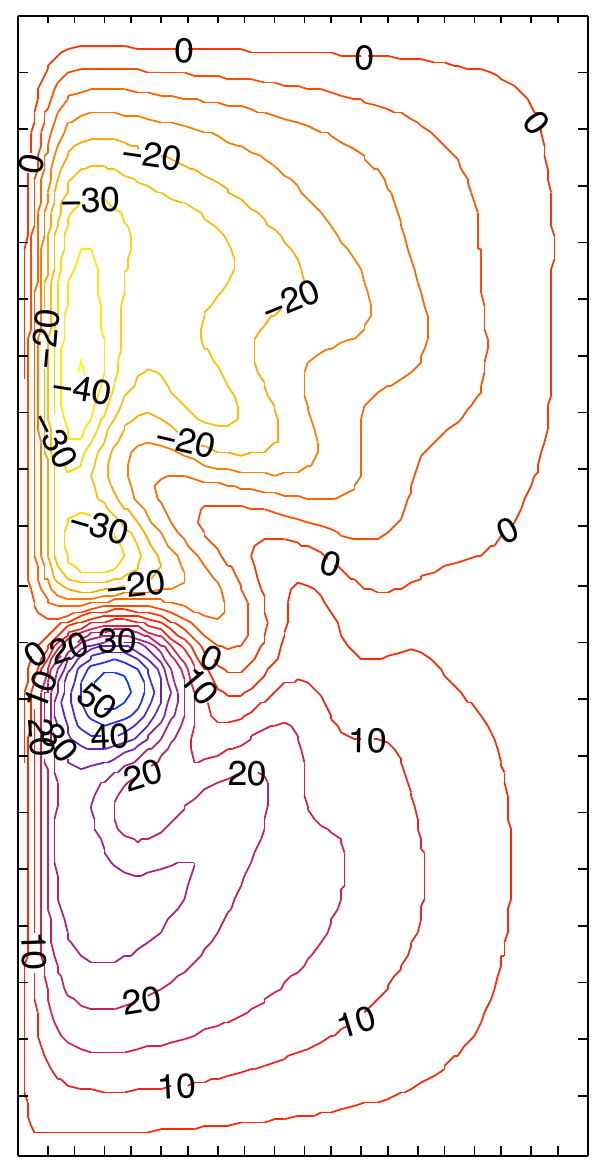}\hspace{7mm}
(c)\includegraphics[width=1.4in]{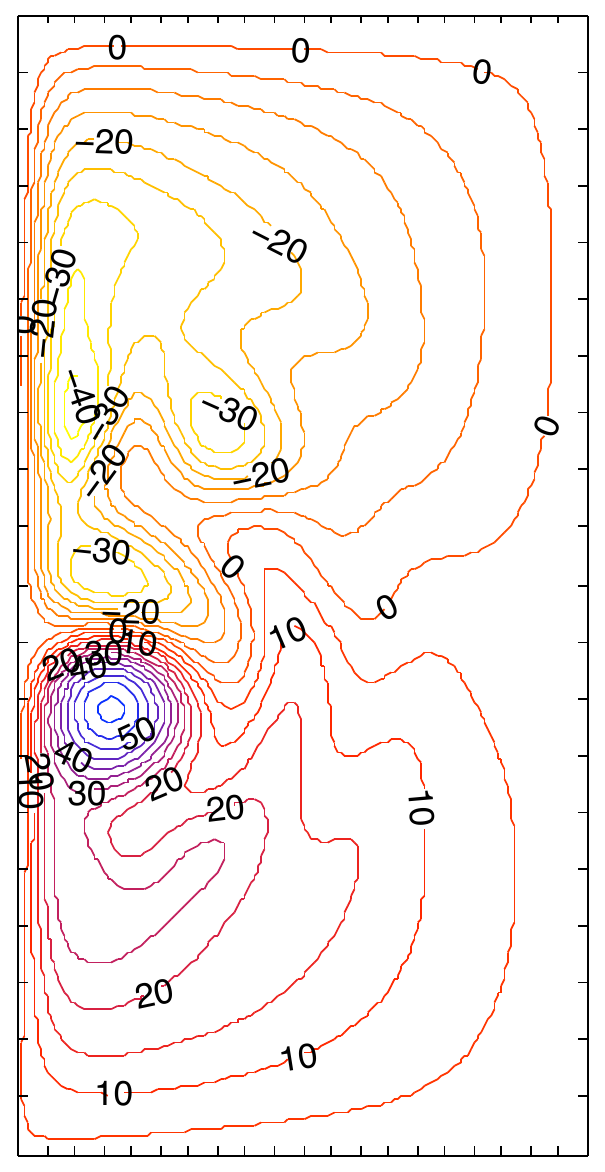}\\ \hspace{7mm}
\end{center}
\caption{Grid refinement study for the proposed formulation and the fractional-step algorithm. The height anomaly of double-gyre model at $t=16$ years. The grid resolutions, from (a) to (c), are $\Delta x=40$ km, $20$ km, and $10$ km, respectively}
\label{fig:f2}
\end{figure}

Figure \ref{fig:f2} shows the height anomaly of the double-gyre model at $t=16$ years calculated by using the new formulation and the fractional-step algorithm. The grid resolutions, from (a) to (c), are $\Delta x=40$ km, $20$ km, and $10$ km, respectively. The dynamics of Figure  \ref{fig:f2}(a) looks different from that of (b) or (c). This is because dynamically the important length scales are only marginally resolved for (a). The length scales are dominated by the first Rossby deformation radius. The Rossby deformation radius for this choice of parameter varies between $45 \sim 80$, from the definition
$$
L_{D} = \frac{1}{f}\left(g_rH_{0}\right)^{1/2}\approx 45 \sim 80 \,\,\text{km}.
$$
The dynamics of the model are dominated by Rossby waves with wave number defined by $\kappa_{R}=1/L_{D}$. If we require that the smallest waves are resolved by the grid spacing $\Delta x$, we must have $\kappa_{R} = \frac{1}{2\Delta x}$ \cite{PJM96}. Hence in this case, the resolution of the grid size must satisfy  $\Delta x \le 20$ km in order to resolve the Rossby waves.

The dynamics of height anomaly, $h-H_0$, of the double-gyre model quickly settles into a quasi-steady-state solution and exhibits strong western boundary currents, as shown in Figure \ref{fig:f2}.  Figure \ref{fig:f2} closely match with Figure 5(a) with $\Delta x=17$ km reported in \cite{PJM96} and Figure 7 reported in \cite{JJG95}. We note that the grid resolution is set to be $\Delta y=\Delta x$ for all our simulations.

%the dynamics of height anomaly computed by the proposed algorithm closely match those in \cite{JJG95,PJM96}.  

\begin{figure}[htpb]
\begin{center}
\includegraphics[width=1.39in]{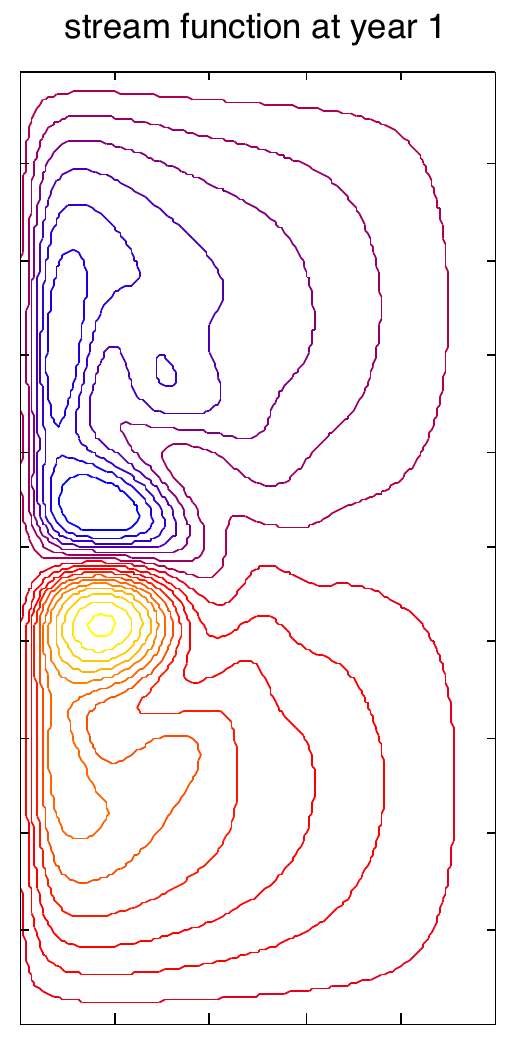}
\includegraphics[width=1.4in]{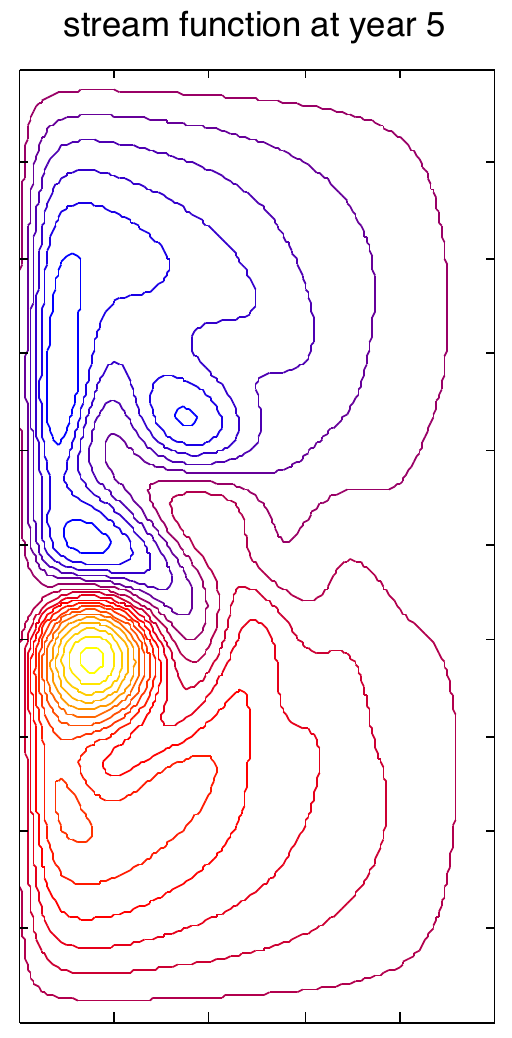}
\includegraphics[width=1.4in]{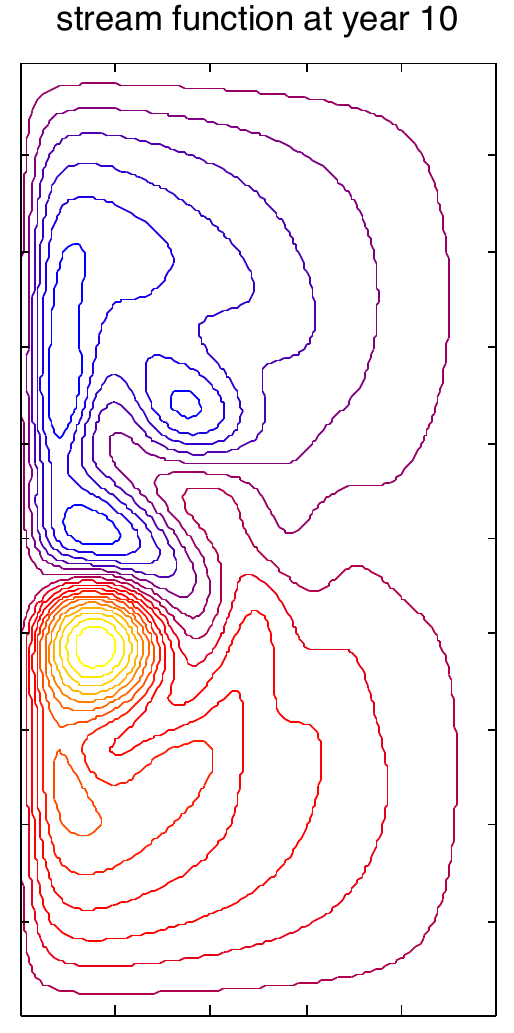}
\includegraphics[width=1.385in]{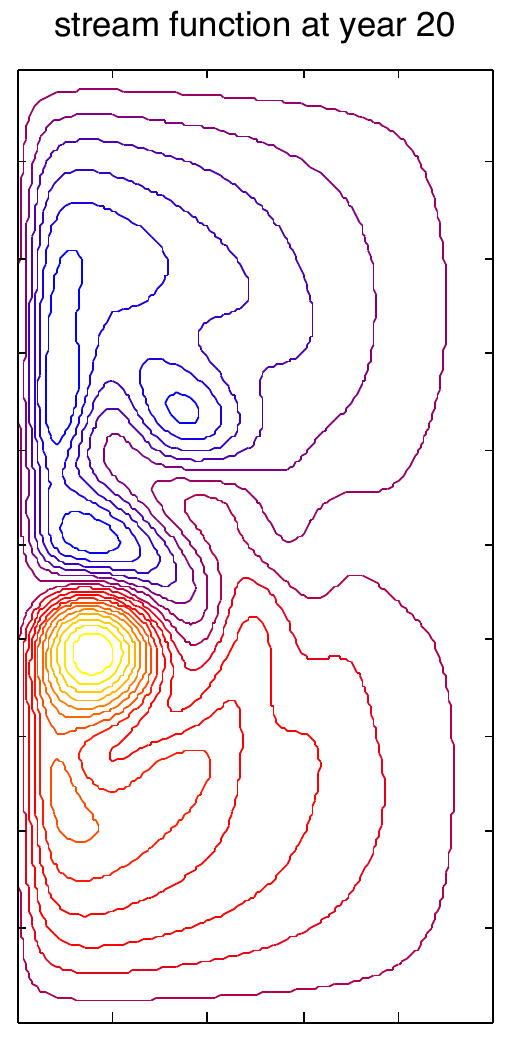}
\end{center}
\caption{From left to right, simulation figures show the contour plots of stream function at year 1, 5, 10, and 20, respectively for a double-gyre model that is under a constant wind forcing described by (\ref{eq:wind_forcing}). The computational domain is $[0, 1000]\times [0, 2000]$ km$^2$. The grid resolution is $\Delta x=\Delta y=10$ km, and the time step is $\Delta t=6$ minutes. For each figure, equally spacing 20 contour lines between $[ -30270, 15944] $ are plotted.}
\label{fig:f3}
\end{figure}

In addition to the height anomaly, we also monitor the velocity field. From left to right, Figure \ref{fig:f3} shows the contour plots of stream function at year 1, 5, 10, and 20, respectively, for the double-gyre model. The wind forcing is described by (\ref{eq:wind_forcing}). The computational domain is again $[0, 1000]\times [0, 2000]$ km$^2$. The grid resolution is $\Delta x=10$ km, and the time step is $\Delta t=6$ minutes. For each simulation figure, 20 contour lines are plotted. We note that the stream-line structures show little difference after year 5 (including year 5). The stream-line structures for year 10 and 20 are almost identical, which provides evidence that the velocity field has reached a quasi-steady-state solution.

Finally, we implement a basic MPDATA algorithm described in \cite{hayder06} for the double-gyre model. Figure \ref{fig:f4} is the comparison of the height anomaly of the double-gyre model after 365 days between the proposed algorithm and the MPDATA implementation. While the structures of the two contour plots are similar, we see that the result from the MPDATA algorithm is more diffusive, even with a mesh that is four-times finer than that for the proposed algorithm.

\begin{figure}[htpb]
\begin{center}
(a) \includegraphics[width=1.95in]{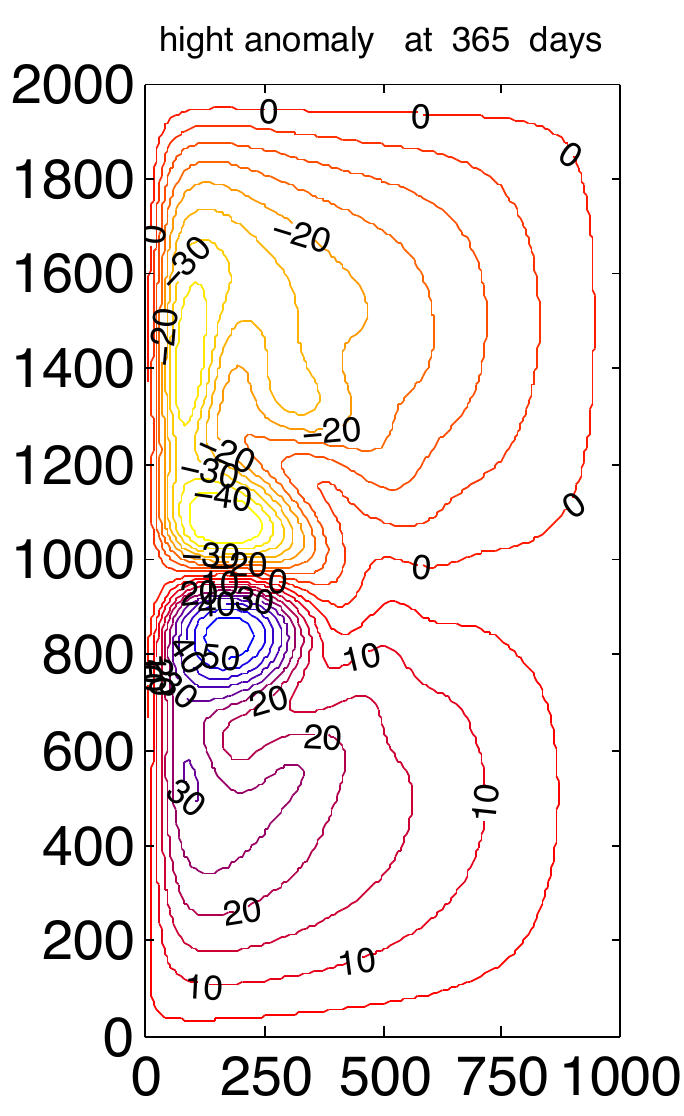}
(b) \includegraphics[width=2in]{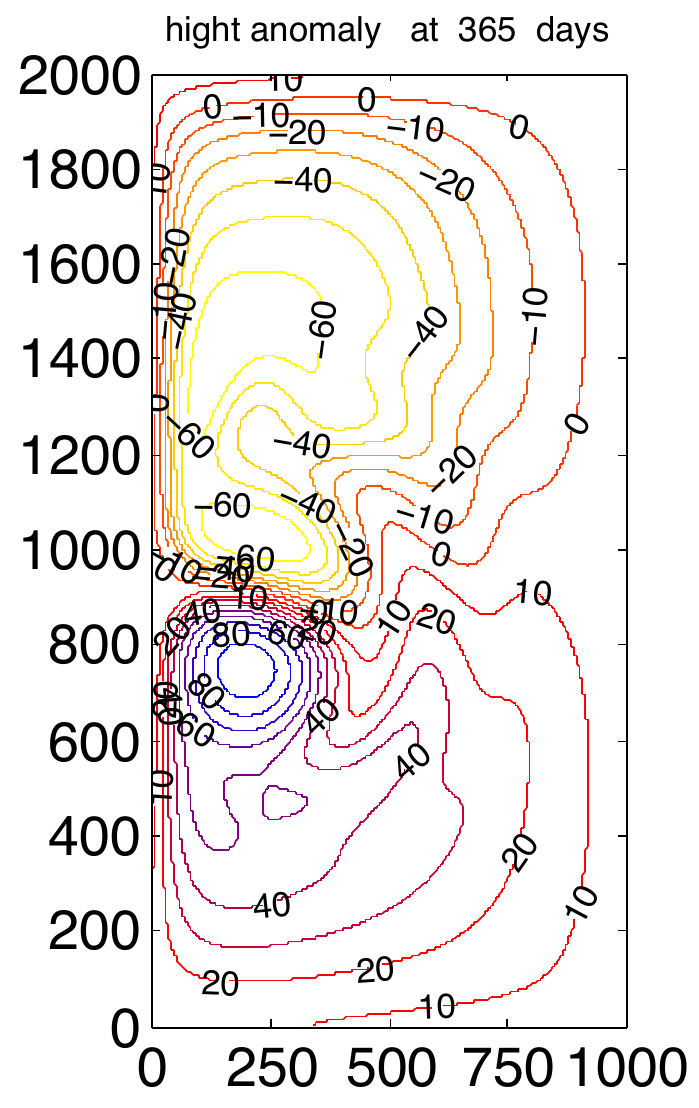}

\end{center}
\caption{The height anomaly of the double-gyre model after 365 days. The enclosed basin with no-slip boundary conditions all around is $[0, 1000]\times [0, 2000]$ $\text{km}^2$. (a) The proposed conservation formulation and the fractional-step algorithm. The grid resolution is $\Delta x= \Delta y =10$ km. $\Delta t= 6$ minutes. (b) Basic MPDATA implementation for the double-gyre model described in \cite{hayder06}. The grid resolution is $\Delta x= \Delta y =2.5$ km. $\Delta t= 0.375$ minutes. }
\label{fig:f4}
\end{figure}

\subsection{Double-gyre model with transport of a pollutant}

The double-gyre shallow-water model has been used as an underlying ocean model for data assimilation \cite{hayder06}. In this section, we use this model to study the circulation of a substance that  initially is randomly distributed in certain areas of a closed ocean basin. This problem is related to transport of pollutant in the ocean, and was previously studied by Xu and Shu \cite{shu06}, using only the hyperbolic shallow-water equations.  To study this problem, we couple the double-gyre shallow-water equations (\ref{eq:double_gyre_conv}) with a two-dimensional scalar advection (transport) equation
%two-dimensional advection-diffusion equation. 
\begin{equation}\label{advection}
%C_t+u\cdot\nabla C = \nu_{d} \nabla^2 C, 
C_t+\bs{u}\cdot\nabla C = 0,
\end{equation}
where $C$ is a substance concentration %$\nu_{d}$ is the diffusion coefficient  of the concentration, 
and $\bs{u}=[u, v]^{T}$ is the velocity field of the double-gyre shallow-water equations. The concentration of the substance is advected by the velocity field of the double-grye shallow-water equations, acting like a scalar tracer. The diffusivity for the scalar tracer is assumed to be very small, so that the diffusion effect of the concentration is negligible. 
The transport equation is solved by the high-resolution wave-propagation algorithm developed in \cite{bib:clawpack}. Because the governing equations are solved in two steps, other than augmenting a conservation equation in the hyperbolic shallow-water equations in Problem A, as suggested in \cite{bib:clawpack}, we solve the color equation( \ref{advection}) in its non-conservation form. The cell-centered value of $C$ is advected by the edge value of the velocity calculated by averaging the adjacent cell-centered velocities. 

We consider a closed ocean basin with dimensions $[0, 1000]\times [0, 2000]$ km$^2$. The basin has been under a constant wind forcing ( \ref{eq:wind_forcing}) for 20 years before the substance is present, and is under the same wind forcing after the substance is present. That is, the quasi-steady-state velocity field, shown in Figure \ref{fig:f3}, is used  as the initial velocity field, and the height field, shown Figure \ref{fig:f2} (b), is used as the initial height field for the double-gyre shallow-water equations. The same parameter values in Table \ref{tab:t2} are used to evolve the double-gyre shallow-water equations. Suppose that the initial values of the concentration are Gaussian random numbers $0 <  C(x,y) \le 1$. We distribute the initial concentration in the following way: Consider two circles with the same radius, $r=150$ km. The centers of the circles are at (500 km, 500 km) and (500 km, 1500 km), respectively. We divide the whole domain into $100\times 200$ grid cells, and assign a random number between 0 and 1 to the center of each grid cell inside the two circles.  Figure \ref{fig:f4} shows the transport of a substance in the basin under the quasi-steady-state velocity field. In the top row, from left to right, the simulation figures show the distribution of concentration at day 0, 50, and 100. In the bottom row, from left to right, the simulation figures show the distribution of concentration at day 150, 240, and 360.  Taken as a whole, Figure \ref{fig:f4} shows that the strong western boundary current drives most of the substance to an area near the western bank. The grid resolution for the simulation is $\Delta x=\Delta y = 10$ km, and the time step is 12 minutes.

\begin{figure}[htbp]
\begin{center}
\includegraphics[width=1.5in]{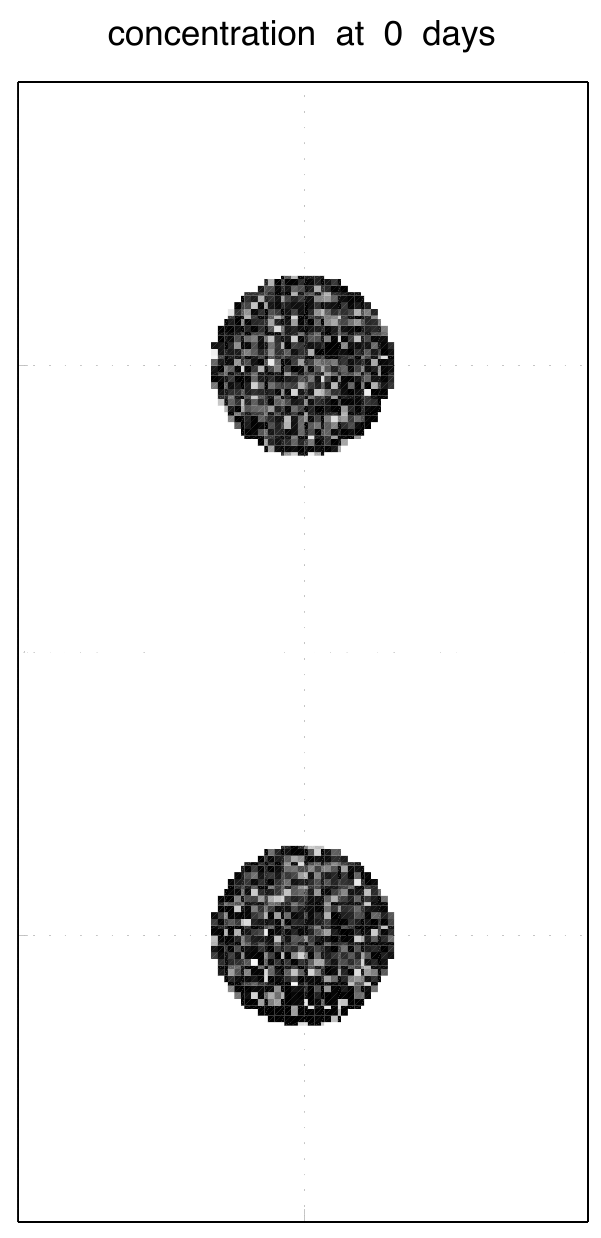}\hspace{2mm}
\includegraphics[width=1.51in]{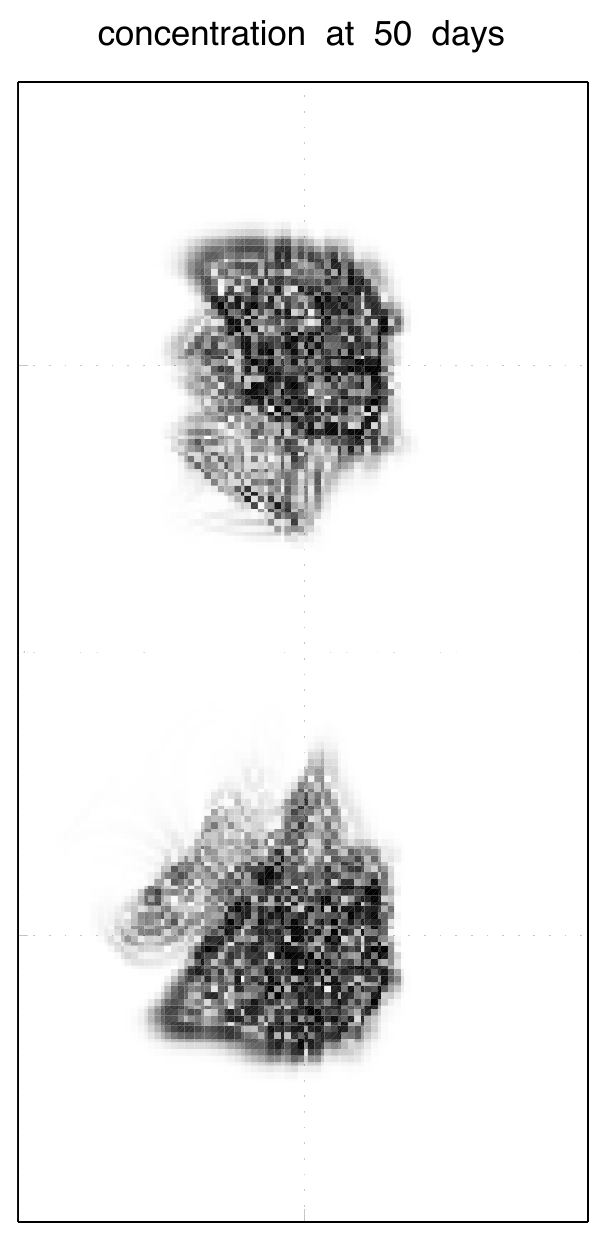}\hspace{2mm}
\includegraphics[width=1.5in]{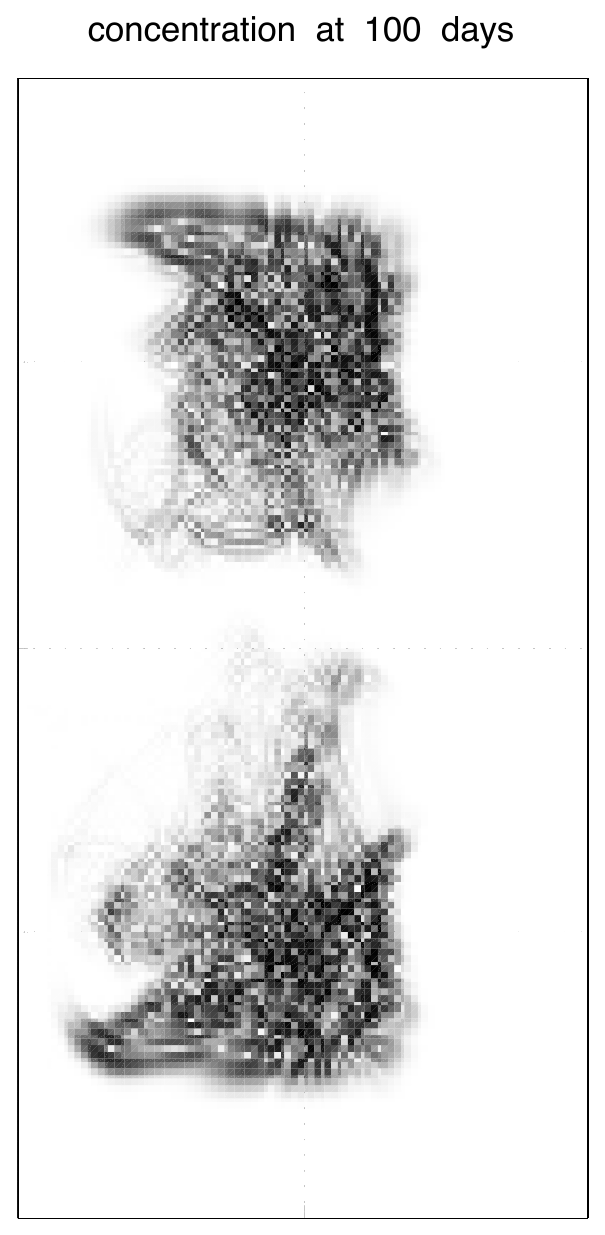}\\%\vspace{10mm}
\includegraphics[width=1.5in]{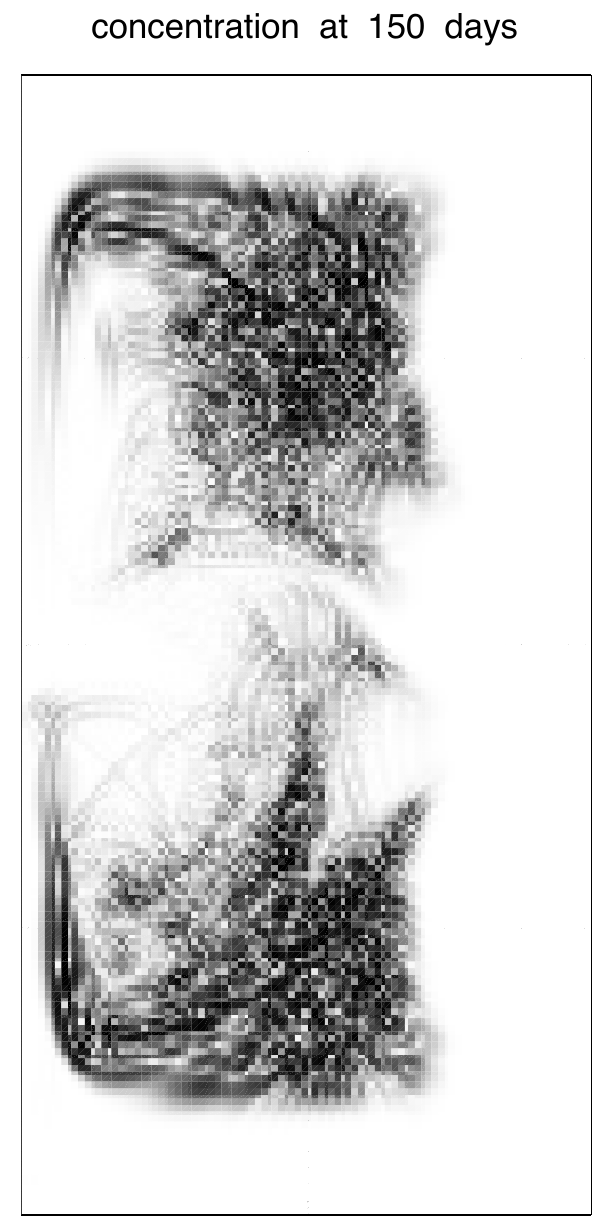}\hspace{2mm}
\includegraphics[width=1.5in]{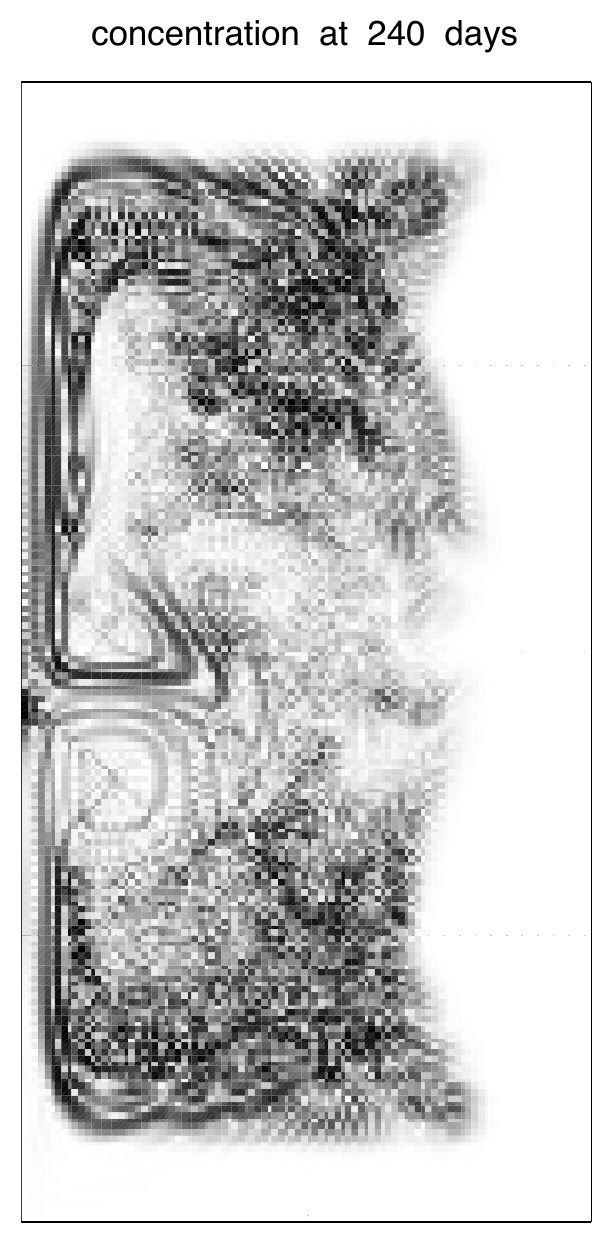}\hspace{2mm}
\includegraphics[width=1.48in]{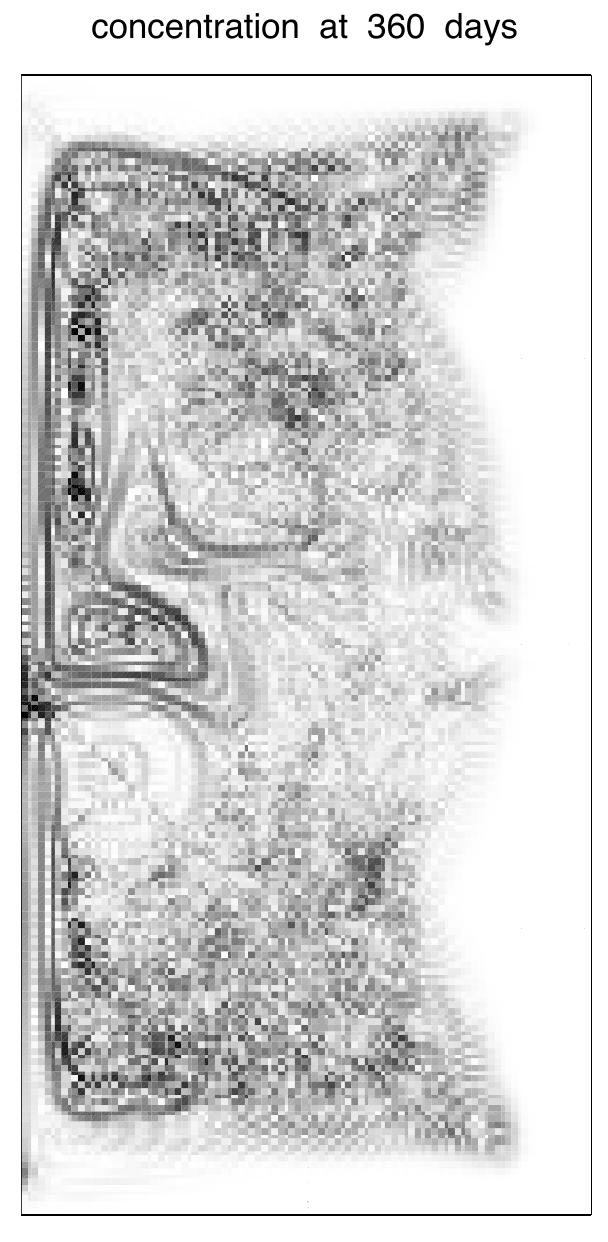}
\end{center}
\vspace{-3mm}
\caption{In the top row, from left to right, the simulation figures show the distribution of a concentration at day 0, 50, and 100. In the bottom row, from left to right, the simulation figures show the distribution of the concentration at day 150, 240, and 360. The concentration is driven by the velocity field of the shallow-water equations. At day 0, the initial concentration is distributed inside the two circles. The concentration values are between zero and one, and are randomly assigned to the center of grid cells inside the circles. Figure \ref{fig:f4} show that the strong western boundary current drives most of the substance to an area near the western bank. The grid resolution for the simulation is $\Delta x=\Delta y = 10$ km, and the time step is 12 minutes. The domain of the basin is $[0, 1000]\times [0, 2000]$ km$^2$. }
\label{fig:f4}
\end{figure}

\section{Conclusion  \label{sec:concluding} }

We present a new formulation for the double-gyre shallow-water model. A fractional-step method is provided to solve the new formulation. The combination of the formulation and the numerical algorithm is proved to be stable and not sensitive to the kinematic viscosity and grid refinement. For traditional methods, stability of the finite difference scheme often depends on the magnitude of  kinematic viscosity. In practice, it is not unusual that to maintain stability, the viscosity needs to be increased as the grid resolution is decreased for those methods \cite{JPM97}. The enslaved finite-difference methods that improves the accuracy for MPDATA could also be sensitive to the viscosity value for certain time integrators when refining meshes. The proposed formulation and the fractional-step method remains stable at a fixed viscosity throughout the gird refinement study. The proposed method is second-order accurate. In the constant wind-forcing example, we demonstrate that the numerical solution converges rather quickly to a quasi-steady-state solution, as long as the Rossby deformation radius is resolved.  Since the high-resolution wave-propagation method that solves the hyperbolic shallow-water equations introduces little numerical dissipation, the proposed fractional-step method is suitable for applications that require small artifical diffusion. Finally, in the last example, we illustrate the flexibility of the proposed method to incorporate other equations for application, such as the transport equation. Especially, when high-resolution is preferable for the monitored quantity in the transport equation. 
%%%% Acknowledgments %%%%%%%%

\section*{Acknowledgments}

The authors thank Zhi (George) Lin for pointing out an error in our earlier numerical implementation.

%%%% Bibliography  %%%%%%%%%%

\end{document}